\documentclass[11pt]{article}

\usepackage{amsfonts}
\usepackage{epic,eepic}
\usepackage{url}
\usepackage{hyperref}
\usepackage{breakurl}

\begin{document}
\bibliographystyle{plain}
\title{{{George Forsythe's last paper
}%
\thanks{\hbox{Presented at the ``Stanford 50'' meeting,
Stanford University, 29 March 2007.}
\hbox{Copyright \copyright 2007--2010, R.~P.~Brent.}
\hspace*{\fill} \hbox{rpb238}}
}}
\author{~\\Richard P.\ Brent\\
	Mathematical Sciences Institute\\
        Australian National University\\
	Canberra, ACT 0200, Australia\\[40pt]
	In memory of\\ 
	George and Sandra Forsythe}
\date{}
\maketitle
\thispagestyle{empty}			%
\vspace*{\fill}

\section*{~}

\begin{quotation}
\noindent
{\em I have a feeling, however, that it is somehow silly to take a random
number and put it elaborately into a power series $\cdots$}
\end{quotation}
\rightline{John von Neumann}
\rightline{Collected Works}
\rightline{Vol.~5, pg.~769}
\bigskip

\vspace*{\fill}
\pagebreak[4]
\section*{Summary}

I will describe von Neumann's elegant idea for sampling from
the exponential distribution, Forsythe's generalization for sampling from a
probability distribution whose density has the form $\exp(-G(x))$, 
where $G(x)$ is easy to
compute (e.g.~a polynomial), and my refinement of these ideas to give an
efficient algorithm for generating pseudo-random numbers with a normal
distribution. Later developments will also be mentioned.

\section*{Background}

Roger Hockney and I are the only people who were lucky enough to have both
George Forsythe and Gene Golub as PhD advisors
(see the {\em  Mathematics Genealogy Project}).  
\smallskip

In my case this came about because Gene went on sabbatical to the UK, and
George took over while he was away. However, I managed to finish before Gene
returned to Stanford.  That was in the days before email, and there was a
mail strike in UK, so communication with Gene was difficult. Perhaps that
helped me to finish quickly, because if Gene had been at Stanford he
probably would have asked me to do more work on the last chapter!
\smallskip

Most of you here today know Gene, but only the older ones will remember
George and Sandra Forsythe, so today I will talk about George Forsythe and
an interesting link back to John von Neumann and the early days of
computers.

\section*{History}

In summer 1949 Forsythe attended some lectures
at UCLA by John von Neumann on the topic of random number generation.
The lectures were part of a Symposium on the 
(then new) 
{\em Monte Carlo method}~\cite[p.~236]{von-Neumann-work}.
It seems that von Neumann never wrote up the lectures, but a fascinating
3-page summary was written by Forsythe and published in 1951.
\smallskip

Forsythe must have continued to think about the topic
because, shortly before he died, he wrote a Stanford report
{\em Von Neumann's comparison method for random sampling from the normal
and other distributions}
(STAN-CS-72-254, dated February 9, 1972).
\pagebreak

This expanded on a brief comment by von Neumann that his
\begin{quotation}
\noindent
{\em method}
{\rm [for the exponential distribution]} 
{\em can be modified to yield 
a distribution satisfying any first-order differential equation.}
\end{quotation}
\rightline{Collected Works {\bf{5}}, 770.}

\section*{Ahrens, Dieter and Knuth}

Forsythe intended that his Stanford report would form the basis of a joint
paper with J.~H.~Ahrens and U.~Dieter, who had discovered related results
independently, and had presented them at Stanford in October 1971. 
\smallskip

After Forsythe died in April 1972, Don Knuth submitted the Stanford report
to {\em Mathematics of Computation}, and it was published with only minor
changes in the October 1972 issue.
\smallskip

This was Forsythe's last published paper, with the possible exception of a
paper by E.~H.~Lee and Forsythe in {\em SIAM Review} (submitted in October
1971 and published in January 1973).
\smallskip

Ahrens and Dieter published a follow-up 
paper in {\em Math.\ Comp.} (1973)~\cite{Ahrens-Dieter73}
and I published an implementation GRAND 
of my improvement of the Forsythe~--~von
Neumann method in {\em Comm.\ ACM} (1974)~\cite{rpb023}. 
\smallskip

That was in the days before {\em TOMS}, when
interesting algorithms were still published in {\em Communications}.

\section*{The problem}

Suppose we want to sample a probability distribution with density 
\[f(x) = e^{-G(x)}\;,\]
where $G(x)$ is some simple function, e.g.~a polynomial.  Von Neumann 
illustrated his idea for the exponential distribution 
\[G(x) = x, \;\; (x \ge 0)\;,\] 
but it also applies to the normal distribution
\[G(x) = x^2/2 + \ln(2\pi)/2\;.\]
The function $f(x)$ satisfies a first-order linear differential equation
\[f' + G'(x)f = 0\;,\]
and conversely. That is why von Neumann made the remark about first-order
differential equations.

\section*{Von Neumann's insight}

The obvious way to generate a sample from the exponential distribution is
to generate a sample $u \in (0, 1]$ from the uniform distribution and then
take
\[ x = -\ln(u)\;.\]
However, the evaluation of $\ln(u)$ is expensive (relative to the cost of
generating $u$ by an efficient uniform random number generator). Also, this
method does not generalize well to the normal distribution, where we would
need to evaluate the inverse of the normal distribution function
\[ \frac{1}{\sqrt{2\pi}} \int_{-\infty}^x \exp(-t^2/2)\,dt\;.\]
Von Neumann's insight 
was that we can generate a random sample using a small number
(on average) of samples from a uniform distribution, and evaluation of
$G(x)$ at a small number of points. There is no need to compute 
any expensive special functions!

\section*{Probability of a run}

Suppose for the moment that $0 \le u_1 = G \le 1$. Generate samples
$u_2, u_3, \ldots$ from the uniform distribution so long as the numbers
are decreasing, and then stop.  In other words, 
find $n \ge 1$ such that
\begin{equation}
G = u_1 > u_2 > u_3 > \cdots > u_n \le u_{n+1}\;.	\label{eq:one}
\end{equation}
The probability that 
\[G > u_2 > u_3 > \cdots > u_n > u_{n+1}\;\;{\rm is}\]
\[\frac{G^n}{n!} =
\frac{{\rm Prob}(\max(u_2, \ldots, u_{n+1}) < G)}{n!}\;,\]
so the probability of~(\ref{eq:one}) is
\[p_n = \frac{G^{n-1}}{(n-1)!} - \frac{G^n}{n!}\;.\]
{\bf Check:} $p_1 + p_2 + \cdots = 1$ by telescoping series, so the 
algorithm terminates with probability~$1$. 
\smallskip

\noindent{\bf Exercise:} The expected value of $n$ is $\exp(G)$.

\section*{The power series for $\exp(-G)$}

What is the probability that our final $n$ is {\em odd}? 
It is just
\[p_1 + p_3 + \cdots = 
1 - G + \frac{G^2}{2!} - \frac{G^3}{3!} + \cdots = \exp(-G)\;.\]
This suggests a {\em rejection method} for 
generating a sample from the distribution with
density $\exp(-G(x))$ on some interval $[a,b]$:
\begin{enumerate}
\item Generate uniform $w \in [a,b]$ and set $u_1 \leftarrow G(w)$.
\item Generate uniform $u_2, u_3, \ldots \in [0, 1]$ until
condition~(\ref{eq:one}) is satisfied $(u_n \le u_{n+1})$.
\item If $n$ is even, return to step 1 (i.e.\ {\em reject} $w$).
\item Return $w$ (i.e.\ {\em accept} $w$).
\end{enumerate}
This works because the probability that $w$ is accepted, i.e.~the
probability that $n$ is odd at step~3, is exactly $\exp(-G(w))$.

\section*{An important condition}

The algorithm only works correctly if $0 \le G(w) \le 1$ on the
interval $w \in [a,b]$.
\smallskip

To apply the idea to the exponential or normal
distributions we have to split the infinite interval $[0,+\infty)$ or
$(-\infty,+\infty)$ into a union of finite intervals $I_k$. Provided the
intervals $I_k$ are small enough, we can use the algorithm to generate
samples from each $I_k$. 
\smallskip

Thus, first select $k$ with the correct probability
\[\int_{I_k} \exp(-G(x))\,dx\;,\]
then use the Forsythe~--~von Neumann algorithm 
to get a sample from $I_k$. 
\smallskip

{\bf Minor detail:} The function $G(x)$ has to be modified by addition 
of a constant to give the appropriate function $G_k(x)$ on the 
interval $I_k$. For example, we could use $(x^2 - a^2)/2$ for the normal
distribution on $[a,b]$, where $0 \le a < b$ and $b^2 - a^2 \le 2$.

\section*{Exponential and normal distributions}

For the exponential distribution, consider the intervals
\[I_k = [(k-1)\ln 2,\; k\ln 2)\;.\]
For convenience on a binary computer,
our sample should lie in $I_k$ with probability $2^{-k}$, $k = 1, 2, \ldots$
We can select $k$
by counting the leading zero bits in a uniform random number (giving $k-1$).
Then we can apply a rejection method to get a sample with the 
correct distribution from $I_k$.  
\smallskip

For the normal distribution, it is convenient to randomly generate the sign,
then consider the interval $[0, \infty)$. We subdivide this interval
into intervals $I_k = [a_{k-1}, a_{k})$ such that $a_0 = 0$ and
\[\sqrt{\frac{2}{\pi}}\int_{a_{k-1}}^{a_{k}} \exp(-x^2/2)\,dx = 2^{-k}\]
for $k = 1, 2, \ldots$, $w$. 
It is easy to precompute a table of the constants
$a_k$. The table is small, since we can neglect probabilities
$2^{-k}$ if $k$ is greater than the wordlength $w$ of the computer.

\section*{Historical notes}

For the exponential distribution,
von Neumann took intervals $I_k = [k-1, k)$ so the probability of sampling
from $I_k$ is $(e-1)/e^k$. He did this because he had a trick for combining
the trials with the selection of intervals.
However, his trick does not generalize to other distributions.
\smallskip

For the normal distribution, Forsythe used intervals defined by $a_0 = 0$
and 
\[a_k = \sqrt{2k-1} \;\;{\rm for}\;\; k \ge 1\;.\]
Presumably he did this because then 
\[a_k^2 - a_{k-1}^2 = 2 \;\;{\rm for}\;\; k \ge 2\;.\]
This choice of $a_k$ 
is what Sandra Forsythe used in her implemention of the algorithm:
\begin{quotation}
\noindent
{\em The correctness of this algorithm $\cdots$ (has) been confirmed in 
unpublished experiments by A.~I.~Forsythe and independently by J.~H.~Ahrens.}
\end{quotation}
\rightline{George Forsythe} %

\section*{Comment}

It is better to use the intervals that I defined, as used in GRAND,
because then we do not need to store a table of probabilities (they are 
just negative powers of~$2$). With my choice it can be shown that
\[a_k^2 - a_{k-1}^2 < 2\ln 2 < 1.39 \;\;{\rm for}\;\; k \ge 1\;.\]
As well as reducing the table size, my choice reduces the expected number
of calls to the uniform random number generator.

\section*{Refinements}

The algorithms proposed by Forsythe and von Neumann were inefficient in the
sense that they used more uniform samples than necessary to generate one
sample from the exponential or normal distribution.
\smallskip
 
The algorithm 
implemented by Sandra Forsythe requires 
(on average) $4.04$ uniform samples per normal sample.
For von Neumann's algorithm the corresponding constant is $5.88$.
\smallskip

Ahrens and Dieter (1973) reduced the constant $4.04$ to $2.54$ (and
even further at the expense of larger tables and more complications).
\smallskip

In my 1974 paper describing GRAND
I showed how $4.04$ could be reduced to $1.38$ 
by using a better subdivision of the infinite interval $[0, +\infty)$
and by not wasting random bits. 
For example, after step 2, 
\[\frac{u_{n+1} - u_n}{1 - u_n}\]
is uniformly distributed and can be used later.

\section*{Further refinements}

In principle, by using larger tables, it is possible to 
reduce the constant to $1 + \varepsilon$ for any $\varepsilon > 0$,
but this would not necessarily give a faster algorithm.
In practice $1.38$ is small enough.

\section*{Later developments}

The idea of {\em rejection methods} was developed by many people to give 
efficient algorithms for sampling from a great variety of distributions~--
see for example the books by Devroye %
and Knuth (Vol.~2).
\smallskip

Specifically for the normal distribution, Forsythe's method (as improved and
implemented in GRAND) is much faster than earlier methods, such as the
Box-Muller and Polar methods.
\smallskip

There are now many different algorithms for the normal distribution,
but I think it is fair to say that none are {\em clearly better} than GRAND.
\smallskip

The differences between the best algorithms are small~-- often there is a
tradeoff between space and time, and the relative speeds depend on the
machine architecture as well as on the choice of uniform random number
generator.

\section*{Wallace's method}

The only method that is clearly much faster than GRAND is Wallace's method,
proposed in 1994 by Chris Wallace. It does not use a uniform random number
generator.  Instead, a pool of normally distributed numbers is maintained
and refreshed by performing orthogonal transformations.
\smallskip

The key observation is that, if $x$ is a vector of $n$
independent, normally distributed numbers, 
then the probability density of $x$,
\[(2\pi)^{-n/2}\exp\left(-(x_1^2 + \cdots + x_n^2)/2\right)\;,\]
is a function of $||x||_2$. i.e.~the distribution has spherical symmetry.
If follows that, if $Q$ is an $n \times n$ orthogonal matrix, then 
\[y = Qx\] 
is another vector of normally distributed numbers,
because $||y||_2 = ||x||_2$.
\smallskip

Wallace's method is interesting and fast, but suffers from some statistical
problems: see my paper in the Wallace memorial~\cite{rpb213}.

\vspace*{\fill}
\pagebreak[4]

\end{document}